\documentclass[11pt]{amsart}
\usepackage{amssymb}

\newtheorem{prop}{Proposition}
\newtheorem{lem}{Lemma}
\newtheorem{thm}{Theorem}

\newcommand{\proj}{\mathbf P}
\newcommand{\grass}{\mathbf G}

\newcommand{\rarr}{\rightarrow}
\newcommand{\oh}{{\mathcal{O}}}
\newcommand{\f}{{\mathcal{F}}}
\newcommand{\com}{\mathbb{C}}

\newcommand{\Z}{\mathbb{Z}}

\newcommand{\Pone}{{\proj^1}}
\newcommand{\F}{\mathbf F}
\newcommand{\R}{\mathcal{HQ}_\bd}
\newcommand{\Rg}{\mathcal{HQ}_d(\grass^r(n))}
\newcommand{\pflag}{\F(n;\bs)}
\newcommand{\bs}{\mathbf{s}}
\newcommand{\s}{(s_1,\ldots,s_l)}
\newcommand{\Rpf}{\R(\pflag)}
\newcommand{\Rcf}{\R(\F(n))}
\newcommand{\bd}{\mathbf{d}}
\newcommand{\dd}{{d_1,\ldots,d_l}}
\newcommand{\hrarr}{\hookrightarrow}
\newcommand{\A}{\mathcal{A}}
\newcommand{\B}{\mathcal{B}}
\newcommand{\ta}{\alpha}
\newcommand{\tb}{\beta}
\newcommand{\ve}{\varepsilon}
\newcommand{\poinc}{\mathcal{P}}
\newcommand{\rrarr}{\twoheadrightarrow}

\newcommand{\eqq}{\stackrel{\sim}{=}}
\newcommand{\sumo}{\oplus}

\newcommand{\bpf}{\noindent {\em Proof.} }
\newcommand{\epf}{\qed \vspace{+10pt}}

\newcommand{\mor}{\mathrm{Mor}}
\newcommand{\Hom}{\mathrm{Hom}}

\begin{document}
\date{\today}
\title{Poincar\'{e} Polynomials of 
Hyperquot Schemes}
\author{Linda Chen}
\maketitle

\section{Introduction}

In this paper, we find generating functions for the Poincar\'{e} polynomials 
of hyperquot schemes for all partial flag varieties.  
These generating functions  
give the Betti numbers of hyperquot schemes, and thus give 
dimension information for
the cohomology ring of every hyperquot scheme.  

Let $\pflag$ denote the partial flag variety corresponding
to flags of the form:
$$V_1 \subset V_2 \subset ... \subset V_l \subset V=\com^n$$
with dim $V_i=s_i$.  The space $\mor_\bd(\Pone, \pflag)$ of morphisms
from $\Pone$ to $\pflag$
 of multidegree $\bd = (d_1,...,d_l)$ can be viewed as the 
space of successive quotients of $V_\Pone$
of vector bundles of rank $r_i$ and degree $d_i$, 
where $r_i:=n-s_i$ and $V_\Pone:=V\otimes \oh_\Pone$ is a trivial
rank $n$ vector bundle over $\Pone$.
  Its compactification,
the hyperquot scheme which we denote $\R=\Rpf$, parametrizes
flat families of successive quotient sheaves  of $V_\Pone$  of 
rank $r_i$ and degree $d_i$.
It is a generalization of Grothendieck's Quot scheme \cite{g}.

There has been much interest in compactifications of moduli spaces of
maps, for example the stable maps of Kontsevich.  
Hyperquot schemes are another natural such compactification.  Indeed, most of
what is known so far about the quantum cohomology of Grassmanians
and flag varieties has been obtained by using Quot scheme compactifications.
They have been used by Bertram 
to study Gromov-Witten invariants and a quantum Schubert calculus \cite{b1}
\cite{b2}  for Grassmannians, and by Ciocan-Fontanine and Kim
 to study Gromov-Witten invariants and the quantum cohomology ring 
of flag varieties and partial flag varieties \cite{k} \cite{cf1} 
\cite{cf2} \cite{fgp}, see also \cite{c1} \cite{c2}.  

The paper is organized in the following way:

In section \ref{quot}, we give some properties of hyperquot schemes, including
a description of the Zariski tangent space to $\R$ at a point.

In section \ref{torus}, we consider a torus action on the hyperquot scheme.  
By the theorems of Bialynicki-Birula,
the fixed points of this action 
give a cell decomposition of $\R$ \cite{bb1}\cite{bb2}. 
The fixed point data is organized to
give a generating function for the topological Euler characteristic of $\R$.

In section \ref{betti}, we use the Zariski tangent space to $\R$ at a point 
as described in section \ref{quot}
to compute tangent weights at the fixed points.
This gives an implicit formula for the Betti numbers of $\R$.
The torus action and techniques are similar to those used by Str{\o}mme in the 
case $l=1$, where $\F(n:s)$ is the Grassmannian $\grass_{s}(n)$ and 
$\R$ is the ordinary Quot scheme \cite{s}.

In section \ref{poincare}, we reorganize the implicit formula for
the Betti numbers in a way that reduces the problem to a purely combinatorial
one.  In particular, we  collect the information into the form of 
a generating function. 

Let $\poinc({\mathbf X}) = \sum_M b_{2M}({\mathbf X})z^M$ denote the 
Poincar\'{e} polynomial of a space ${\mathbf X}$.  It is classically known
that $\poinc(\pflag)$ is equal to the following generating
function for the Betti numbers of the partial flag variety:
  $$\poinc(\pflag) =\sum_M b_{2M}(\F)z^M = \frac{\prod_{i=1}^n (1-z^i)}
    {\prod_{j=1}^{l+1} \prod_{i=1}^{s_j-s_{j-1}} (1-z^i) }$$
with $s_{l+1}:=n$ and $s_0:=0$.

Defining $f^{i,j}_k:=1-t_i\cdots t_jz^k$, the main result is:

\begin{thm}
\label{Rpf}
\begin{eqnarray*}
\lefteqn{\sum_\dd \poinc(\Rpf)t_1^{d_1}...t_l^{d_l}= }\\
& &  \poinc(\pflag) \cdot
\prod_{1\leq i \leq j \leq l}
\prod_{s_{i-1}<k\leq s_i}
\left(\frac{1}{f^{i,j}_{s_j-k}}\right)
\left(\frac{1}{f^{i,j}_{s_{j+1}-k+1}}\right)
\end{eqnarray*}

\end{thm}

In section \ref{special}, we  discuss the special cases of the ordinary
Quot scheme and of the hyperquot scheme for complete flags, and provide some specific
examples.

\section{Hyperquot Schemes}
\label{quot}

We fix some notation. 
Let $V=\com^n$ be a complex n-dimensional vector space. 
Let $\F:= \pflag$ denote the partial flag variety corresponding
to flags of the form:
$$V_1 \subset V_2 \subset ... \subset V_l \subset V$$
where $V_i$ is a complex subspace of dimension $s_i$.
We have $s_0:=0 < s_1<...<s_l< s_{l+1}:=n$.  Define $r_i:=n-s_i$.
As a special case, let $\F(n)=\F(n;1,2,...,n-1)$ 
denote the complete flag variety, with $\Rcf$ the
corresponding hyperquot scheme for complete flags.  Also note that the Grassmannian 
parametrizing $r$-dimensional quotients of $V$, is also
a special case, $\grass^r(n)=\F(n;n-r)$.

For any space $T$, let $V_T$ denote the trivial rank $n$ vector bundle on $T$, 
i.e. $V_T := V \otimes \oh_T$ .

Consider a functor $\f_\bd$ from the category of schemes to the 
category of sets.  For a scheme $T$, 
$\f_\bd (T)$ is defined to be the set of flagged quotient sheaves
$$V_{\Pone \times T}\rrarr \B_1 \rrarr \cdots \rrarr \B_l$$
with each $\B_i$ flat over $T$ with Hilbert polynomial
$\chi(\Pone,(\B_i)_t(m)) = (m+1)r_i+d_i$ on the fibers
of $\pi_T: \Pone \times T \rightarrow T$.  This last condition
requires that $\B_i$ be of rank $r_i$ and relative degree $d_i$ over $T$,
so  that  for any $t\in T$, $(\B_i)_t$ is of degree $d_i$.

It is proven that the functor $\f_\bd$ is represented by the projective
scheme $\R=\Rpf$ following the ideas of Grothendieck and
Mumford \cite{cf2}\cite{g}\cite{m}.  It has also been 
described in a different way by Kim \cite{k}, 
as a closed subscheme of a product of 
Quot schemes.  Kim also proves the following result:

\begin{thm}[Kim]
\label{Rpfproj}
$\Rpf$ is an irreducible, rational, nonsingular, projective variety
of dimension 
$$\sum_{i=1}^l d_i(s_{i+1}-s_{i-1}) + \dim(\pflag)$$
 \end{thm}
\noindent with
$s_0=0,s_{l+1}=n$ and $ \dim(\pflag) =\sum_{i=1}^l (s_{i+1}-s_i)s_i.$
In particular, the theorem states that
 $\dim \Rcf = 2|d|+ \binom{n}{2}$ where 
$|d| = \sum_{i=1}^l d_i.$

Associated to  $\Rpf$ is a universal sequence of sheaves on 
$\Pone \times \R$ of
successive quotients of sheaves, each of which is flat over $\R$.
$$ V_{\Pone\times \R} 
\twoheadrightarrow B_1 \twoheadrightarrow \cdots \twoheadrightarrow B_l.$$

Define $A_i$ as the kernel of $V_{\Pone\times \R} \rarr B_i$. Each
$A_i$ is flat over $\R$, and it
is an easy consequence of flatness that each $A_i$ is locally free.
Thus, we have the following universal sequence on  $\Pone \times \R$:
\begin{equation}
\label{univ} 
A_1 \hrarr A_2 \hrarr \cdots \hrarr A_l \hrarr V_{\Pone\times \R} 
\twoheadrightarrow B_1 \twoheadrightarrow \cdots \twoheadrightarrow B_l.
\end{equation}
with $A_i$ of rank $s_i$, $B_i$ of rank $n- s_i$.  Denote the
inclusion maps by $\gamma_i: A_i \hrarr A_{i+1}$ and the surjections
by $\pi_i:B_{i-1} \rarr B_i$ for each $1\leq i\leq l$.
Here, we define  $A_{l+1}=B_0= V _{\Pone\times \R}$ and
$A_0=B_{l+1}=0$.  
The map  $\gamma_i: A_i\hrarr A_{i+1}$ is an inclusion of sheaves, 
not an inclusion of bundles.

The following proposition, proved by Ciocan-Fontanine following the
ideas in Kollar's work on  Hilbert schemes,
 determines the Zariski tangent space of $\R$ at a point  \cite{ko}.

\begin{prop}
\label{TR}
Let $x\in \R$ correspond to successive quotients and subsheaves of $V_\Pone$:
$$\A_1 \hrarr  \cdots \hrarr \A_l \hrarr V_{\Pone\times \R} 
\twoheadrightarrow \B_1 \twoheadrightarrow \cdots \twoheadrightarrow \B_l.$$
Then we have the following exact sequence:
$$
0 \rarr (T_{\R})_x \rarr \bigoplus_{i=1}^l \Hom(\A_i,\B_i) \stackrel{d}{\rarr}
\bigoplus_{i=1}^{l-1} \Hom(\A_i,\B_{i+1})\rarr 0$$
where $(T_{\R})_x$ is the Zariski tangent space to $\R$ at the point $x$,
and $d$ is the restriction of the difference map given by $d(\{\phi_i\})
= \{\pi_{i+1}\circ\phi_i - \phi_{i+1}\circ\gamma_i\}$.
\end{prop}

\section{A Torus Action}
\label{torus}

In this section, we use the torus action introduced by Str{\o}mme in Theorem
3.6 of \cite{s}.
In the case of the ordinary Quot scheme, we obtain the same description of
the fixed points as Str{\o}mme, but with slightly different notation.  Our
description allows us to provide a full description of the fixed point locus
of the hyperquot scheme under this torus action.

Consider a maximal (n-dimensional) torus $T$ in $GL(V)$ 
which acts on $V$ and hence
induces an action on subsheaves of $V_\Pone$.
As a $\com T$-module, $V$ splits as a direct sum of one-dimensional subspaces,
$\sumo_{i=1}^n W_i$.  Denote $\oh_i := W_i \otimes \oh_\Pone$.

For $f_k \in H^0(\Pone,\oh_\Pone(d_k))$, a form of degree $d_k$, let
$ f_k \oh_i(-d_k) \hrarr \oh_i$
denote the sheaf $\oh_i(-d_k)$ defined by the section $f_k$.

\begin{lem}
A locally free subsheaf $\mathcal{S}\hookrightarrow V_\Pone$ of rank $s$
and degree $-d$ is fixed by the
action of $T$ if and only if it is of the form
$$\mathcal{S} = \bigoplus_{k=1}^{s} f_k \oh_{c_k}(-d_k)$$
where $d_k$ are nonnegative integers such that $\sum_k d_k = d$, 
$1\leq c_1<  \cdots <c_s\leq n$, and $f_k$ is a homogeneous form in $X$ and $Y$
of degree $d_k$.

\end{lem}

\bpf
Since $T$ acts
with different weights on each $\oh_i$,  $\mathcal{S}\hookrightarrow V_\Pone$
is a fixed point of  $T$ if and only if
$\mathcal{S} = \bigoplus_{i=1}^n \mathcal{S}_i$
 where $\mathcal{S}_i:=\mathcal{S} \cap \oh_i$.
Since $\mbox{rank}(\mathcal{S}) = s$, $\mathcal{S}_i \neq \emptyset$ 
for exactly $s$
such $i$, say  $\mathcal{S}_{c_k} \neq \emptyset$ for a sequence 
$1\leq c_1<  \cdots <c_s\leq n$.

Since we have $\mathcal{S}_{c_k}\hrarr \oh_{c_k}$, we know that we can write
$\mathcal{S}_{c_k}=\oh_{c_k}(-d_k)$ for some nonnegative integer $d_k$.  

Threfore we  have:

$$\mathcal{S}=\bigoplus_{k=1}^s \mathcal{S}_{c_k} = 
v\bigoplus_{k=1}^s \oh_{c_k}(-d_k),$$ where 
deg $\mathcal{S}_{c_k} = -d_k$ with $\sum_{k=1}^s d_k = d$.  
Since the  inclusion of $\mathcal{S}_{c_k}$ is given by some section $f_k$,
which is a polynomial of degree $d_k$ in $X$ and $Y$, the lemma is proven.

Let $T'$ be the one-dimensional torus which acts on
$H^0(\oh_\Pone(1))$ by $X \mapsto tX$ and $Y \mapsto t^{-1}Y$.  Then $T'$ acts
on $\Pone$ and hence on subsheaves of $V_\Pone$.  
Thus under the action of the product torus $T\times T'$,
the  fixed points are subsheaves 
$\mathcal{S}\hrarr V_\Pone$ fixed by both $T$ and $T'$.
A point $(\oplus_{k=1}^s f_k \oh_{c_k}(-d_k)\hrarr V_\Pone)$
is fixed by $T'$ if and only if each $f_k$ is a monomial in $X$ and $Y$.

Therefore, we have proven:
\begin{lem}
\label{gfix}
A locally free subsheaf $S\hrarr V_\Pone$ of rank $s$ and degree $-d$
is fixed under the action of ${T \times T'}$  
if and only if it is of the form
$$\bigoplus_{k=1}^s X^{a_k} Y^{b_k} \oh_{c_k}(-a_k-b_k).$$
\end{lem}
\noindent
Here, $\oh_i$ denotes the $i$th component of the trivial rank $n$ vector
bundle $V_\Pone$, and $(a,b,c)$ are sequences of $s$
nonnegative integers satisfying:

\begin{enumerate}
\item
$\sum_{k=1}^s a_k + b_k = d$

\item $1\leq c_1<\cdots < c_s\leq n$

\end{enumerate}

{\noindent {\em Remark.} }
This combinatorial data is equivalent to the fixed 
point data of Str{\o}mme.  For an  element $(\alpha,\beta,\delta)$ as in \cite{s},
let $\delta_{c_1}=...\delta_{c_{n-r}} = 1$
be the nonzero elements, with $1\leq c_1 <... < c_{n-r} \leq n$.
Then the sequence $(\alpha,\beta,\delta)$ corresponds to the sequence 
$(a,b,\mathbf{c})$.

\subsection{A torus action on $\R$}
\label{torusR}
Note that
a point of $\R$ can be given by successive subsheaves over $\Pone$ of 
$V_\Pone = \oplus_{i=1}^{n} \oh_i$.

Let $T$ and $T'$ be as above.
Since the actions of $T$ and $T'$ extend to  actions
on $\Rpf$, we have an action of $T\times T'$ on $\R$.

Using the same methods as used in section \ref{torus}, 
we find the fixed points of $\R$ under this action. 

A point of $\R$ can be  given by a sequence of subsheaves
$$ \A_1 \hrarr \cdots \hrarr \A_l \hrarr V_\Pone$$
where rank $\A_i = s_i$ and $\deg \A_i = -d_i$.
Let $\A$ denote this sequence $\{\A_i\}_{i=1}^l$.  Then $\A$ is fixed
by the action of $T\times T'$ if and only if each $\A_i\hrarr V_\Pone$ 
is fixed and the inclusions $ \A_i \hrarr \A_{i+1}$ hold.  
By Proposition \ref{gfix}, $\A_i$ is fixed when
$$\A_i =\bigoplus_{j=1}^{s_i} \A_{i,j}:=
\bigoplus_{j=1}^{s_i} X^{a_{i,j}} Y^{b_{i,j}} 
\oh_{c_{i,j}}(-a_{i,j}-b_{i,j})$$
where $\sum_{1\leq j\leq s_i} a_{i,j}+b_{i,j} = d_i$ and 
$1\leq c_{i,s_{i-1}+1}< \cdots <c_{i,s_i} \leq n$.  Here, we denote
by $\oh(-a-b)$ the line bundle on $\Pone$ given by global section 
$X^aY^b$.

The inclusion $ \A_i \hrarr \A_{i+1}$ holds under exactly the following
conditions:
\begin{enumerate}
\item
$c_{i,j} = c_{i+1,j}$ whenever $1\leq j\leq s_i$.

\item  \begin{equation}
\label{subsheaves}
\oh_{c_{i,j}}(-a_{i,j}-b_{i,j}) \hrarr  
\oh_{c_{i,j}}(-a_{i+1,j}-b_{i+1,j})
\end{equation}
is an inclusion of sheaves.
\end{enumerate}

Let $S:=S(n;s_1,...s_l)$ be the subset of $S_n$ consisting of  permutations 
$\sigma\in S_n$  such that if
$\sigma(i)<\sigma(i+1)$ unless $i\in (s_1,...s_l)$.  
More explicitly, an element $\sigma\in S$ is such that 
\begin{equation}
\label{sigma}
\sigma(s_{i-1}+1)<\cdots <\sigma (s_i) \mbox{ for } 1\leq i\leq l.
\end{equation}
Therefore every sequence $\{c_{i,j}\}$ corresponds to an element
$\sigma\in S$ by $c_{i,j}=\sigma(j)$, and this correspondence is
a bijection by the first condition above.

The second condition  gives  conditions on the sequences of
nonnegative integers $a$ and $b$.
The inclusion of sheaves $\oh(-a_{i,j}-b_{i,j}) \hrarr \oh$ is given
by the global section $X^{a_{i,j}}Y^{b_{i,j}}$ and
$\oh(-a_{i+1,j}-b_{i+1,j}) \hrarr \oh$ is given by the global section
$X^{a_{i+1,j}}Y^{b_{i+1,j}}$.  Therefore  
(\ref{subsheaves}) gives an inclusion of subsheaves if and only if
that inclusion is given by global sections 
$$X^{a_{i,j}-a_{i+1,j}}Y^{b_{i,j}-b_{i+1,j}}$$ so that the conditions
$a_{i,j}-a_{i+1,j}\geq 0$ and $b_{i,j}-b_{i+1,j}\geq 0$ must hold.

Let $P$ be the set of 
 $(a,b,\sigma)$ such that $a$ and $b$ are sequences of nonnegative integers
$a_{i,j}$ and $b_{i,j}$ with $1 \leq i \leq l, 1\leq j\leq s_i$, and
$\sigma\in S(n;s_1,...,s_l)$ which satisfy:

\begin{enumerate}
\item
$a_{i,j}\geq a_{i+1,j}$

\item
$b_{i,j}\geq b_{i+1,j}$

\item
$\sum_{j=1}^{s_i} a_{i,j} + b_{i,j} = d_i$

\end{enumerate}

For $(a,b,\sigma)\in P$, define
$$\A_{i,j} =  \left\{ \begin{array}{cl} 
\oh_{\sigma(j)}(-a_{i,j}-b_{i,j}) &
\mbox{ for } 1\leq i\leq l\mbox{ and } 1\leq j \leq s_i\\ 0 &
\mbox{ otherwise} \end{array}\right. $$

This defines
sequences of subsheaves $\mathcal{A}_{i,j} \hrarr \mathcal{A}_{i+1,j}$.
Let $\A_i=\oplus_{j=1}^{s_i} \mathcal{A}_{i,j}$.
Let $r(a,b,\sigma)\in \R$ be the associated 
 flag of subsheaves of $V_\Pone$.  We have proven:

\begin{prop}
$r:P\rarr \R$ is a bijection onto $\R^{T\times T'}$.
\end{prop}

Let $\B_{i,j}:=V_\Pone/\A_{i,j}$ be
the corresponding quotient so that we have the following short exact
sequences:
\begin{equation}
\label{ses}
0 \rarr \A_{i,j}\rarr V_\Pone \rarr \B_{i,j} \rarr 0
\end{equation}
Similarly, define the sheaves $\B_i =\oplus_{j=1}^{s_i} \mathcal{B}_{i,j}$.

\subsection{Euler characteristic}

Under a generic choice of one-dimensional 
subtorus $\Gamma\subset T \times T'$, $\R$ has isolated fixed points under
the action of $\Gamma$, with $\R^\Gamma =\R^{T\times T'}$.  
We know by  Theorem \ref{Rpfproj} that  $\Rpf$ is a nonsingular
complex projective variety.  The odd cohomology of $\R$ vanishes since the
fixed point locus $\R^\Gamma$ is finite \cite{bb1} \cite{bb2}. 
Therefore, the Euler characteristic is the number of fixed points.  
The fixed point data
has been collected into the combinatorial data of the set $E$, so that
the Euler characteristic is the cardinality of $E$.
In this section, we  find a generating function for the Euler
characteristics of the hyperquot schemes.  We prove:

\begin{thm}
\label{eulerthm}
$$\sum_{d_1,...d_l} \chi(\Rpf) t_1^{d_1}...t_l^{d_l}= \left(\#(S)\right)
\left(\prod_{1\leq i \leq j \leq l}
\frac{1}{(1-t_i...t_j)^{s_i-s_{i-1}}}
\right)^2
$$
\end{thm}
\noindent
The cardinality of the set of permutations $S=S(n;s_1,\cdots,s_l)$ is:
$$\#(S) = \prod_{i=0}^l \binom{n-s_{i}}{s_{i+1}-s_{i}}$$
with $s_0=0, s_{l+1}=n$.

\bpf
Let $\ta_{i,j}=a_{i,j}-a_{i+1,j}$ and $\tb_{i,j}=b_{i,j}-b_{i-1,j}$,
where we let $a_{l+1,j}=b_{l+1,j}=0$.  Then
$a_{i,j}= \sum_{k=i}^l \ta_{k,j}$ and $b_{i,j}= \sum_{k=i}^l \tb_{k,j}$.
Note that while the variables $(a,b)$
satisfy various inequalities, the nonnegative integers of $(\ta,\tb)$
are independent. Consider the set $P'$ of elements
$(\ta,\tb,\sigma)$ with $\sigma\in S$ satisfying:
\begin{equation}
\label{rel}
\sum_{k=i}^l \sum_{j=1}^{s_i} \ta_{k,j} + \tb_{k,j} = d_i.
\end{equation}

Let $P''$ denote the set of pairs $(\ta,\tb)$ of natural numbers 
satisfying the linear relations (\ref{rel}).
We have constructed a bijection between $E$ and $E'$, so that $\chi(\Rpf)
= \#(P) = \#(P') = \#(S)\#(P'')$

We see that 
\begin{eqnarray*}
\lefteqn{
\left(\prod_{1\leq i \leq j \leq l}
\frac{1}{(1-t_i...t_j)^{s_i-s_{i-1}}}
\right)^2  } \\
&=&\prod_{1\leq i \leq j \leq l} 
\prod_{s_{i-1}<k\leq s_i}
\sum_{\alpha_{k,j}\in\mathbb{N}} (t_i...t_j)^{\alpha_{k,j}}
\sum_{\beta_{k,j}\in\mathbb{N}} (t_i...t_j)^{\beta_{k,j}}\\
 &=& \sum_{\alpha_{k,j}} \sum_{\beta_{k,j}}
\prod_{i=1}^l t_i^{\sum_{k=i}^l \sum_{j=1}^{s_i} \ta_{k,j} + \tb_{k,j}}
\end{eqnarray*}
so that each set of natural numbers $(\ta,\tb)$ satisfying the relations 
(\ref{rel}) contributes exactly one to the coefficient of $t_1^{d_1}\cdots
t_l^{d_l}$.  This proves the theorem.

\section{An implicit formula for the Betti numbers}
\label{betti}

Let $N>>0$ be a large integer.  Let
 $\Gamma\subset T \times T'$ be the one-dimensional subtorus which
acts on $\oh_i$ by $t\cdot v = t^{Ni}v$ and acts on $H^0(\Pone,\oh(1))$ by
$X\mapsto tX, Y\mapsto t^{-1}Y$.  It is this
$\com^*$ action of $\Gamma$ on $\Rpf$ which is used.

In order to find information about the contribution of the various
fixed points to the Betti numbers, we must
consider our $\com^*$ action  of $\Gamma$ on $\R$ more carefully.  

Fix the following notation.  Consider $\Gamma$ a $\com^*$-action 
on a scheme $X$.  If $x\in X$ is a fixed point
of this action, and  $E$ is a bundle over
$X$, denote by $E_x^+$ the $\com\Gamma$-
submodule of $E(x)$ where $\Gamma$ acts with positive weights.
In particular, the theorems of Bialynicki-Birula state that
for isolated fixed points $\{x_i\} \in X_\Gamma$,
there is a cell decomposition of $X$ given by the orbits $X_i$ of $x_i$,
with $\dim_\com X_i = \dim_\com T_X(x_i)^+$ \cite{bb1} \cite{bb2}.

We  first  compute the tangent weights at
the fixed points of the $\com^*$ action of $\Gamma$ on $\R$.
For $\sigma\in S$, define $\epsilon_{i,j}^\sigma$ for $1\leq i,j \leq n$ by
$\epsilon^\sigma_{i,j} = 1$ if $\sigma(i)<\sigma(j)$, 
$\epsilon^\sigma_{i,j} = 0$
otherwise.  
Note that $\epsilon_{i,j}+ \epsilon_{j,i} =1$
for $i\neq j$.  Define $\varepsilon_{i(k,k']}^\sigma$ for $1\leq i,k,k' \leq l$
by $\ve^\sigma_{i(k,k']} = \sum_{s_k<j\leq s_k'} \epsilon_{i,j}$
and similarly define 
$\ve^\sigma_{(k,k']j}=\sum_{s_k<i\leq s_k'} \epsilon_{i,j}$.
If the permutation $\sigma$ is understood, it may be suppressed.

Define a map $h:P\rarr\Z$ by:
\begin{eqnarray}
\label{hpf}
h(a,b,\sigma) &=& \sum_{i=1}^l \sum_{k\leq s_i} 
(a_{i,k}+b_{i,k}+1)\ve^{\sigma}_{k(i,i+1]} \\
\nonumber &+& 
\sum_{i=1}^l \sum_{k\leq s_i} (a_{i,k}+b_{i,k})\ve^{\sigma}_{(i-1,i]k} + 
\sum_{i=1}^l \sum_{s_i<k\leq s_{i+1}} b_{i,k}
\end{eqnarray}

Then we have the following implicit formula for the Betti numbers
of the hyperquot scheme:
\begin{prop}
\label{bpf}
For any (nonnegative) integer $M$,
$$b_{2M}(\Rpf) =\mbox{rank } A_M(\R) = \#(h^{-1}(M)).$$
\end{prop}

\bpf
If we show that 
$h(a,b,\mathbf{e}) = \dim_\com (T_{\R} (r(a,b,\sigma))^+)$, then by applying 
the theorems of Bialynicki-Birula, we will be done.  Let $\mathcal{A}$ be
the sequences of subsheaves associated to $r(a,b,\sigma)\in \R$ as in
section \ref{torusR}.
From the short exact sequence on $\R$ given in Lemma \ref{TR} we have:

$$0 \rarr T_{\R}(\A) \rarr \bigoplus_{i=1}^{l} \Hom(\A_i,\B_i) \rarr
\bigoplus_{i=1}^{l-1} \Hom(\A_i,\B_{i+1})\rarr 0.$$

Therefore the tangent weights of $T_{\R}(\A)$  are those obtained by
removing the  weights of the quotient term from those
of the middle term.  In particular, the positive weights are
also be obtained this way.  More precisely, we can say that 
\begin{eqnarray*}
\lefteqn{ \dim_\com(T_{\R}(\A))^+ = }\\
& & \dim_\com\left(\bigoplus_{i=1}^{l} \Hom(\A_i,\B_i)\right)^+ - 
\, \dim_\com\left(\bigoplus_{i=1}^{l-1} \Hom(\A_i,\B_{i+1})\right)^+.
\end{eqnarray*}

Define maps $h_{i,j}^k:P\rarr \Z$ 
for $0\leq i \leq j\leq l, k=1,2,3$ as follows:
\begin{eqnarray}
\label{definehk}
& &h_{i,j}^1(a,b,\sigma) =  \sum_{k\leq s_i}
(a_{i,k}+b_{i,k}+1)\ve^{\sigma}_{k,(j,l+1]}
\\ \nonumber & &h_{i,j}^2(a,b,\sigma) = 
\sum_{k\leq s_j}(a_{j,k}+b_{j,k})\ve^{\sigma}_{(0,i],k} \\
\nonumber & &h_{i,j}^3(a,b,\sigma) = 
\sum_{k\leq s_i} b_{j,k}
\end{eqnarray}
where we allow zero maps when appropriate.

The proposition is then an immediate consequence of the following
two lemmas.

\begin{lem}
\label{kk}
For any $i\leq j$,
  $\dim_\com(\Hom(\A_i,\B_j))^+ = (h_{i,j}^1+h_{i,j}^2+h_{i,j}^3)(a,b,\sigma)$.

\end{lem}

\begin{lem}
\label{sum}
$$h(a,b,\sigma) = \sum_{i=1}^l (h_{i,i}^1+h_{i,i}^2+h_{i,i}^3)(a,b,\sigma)-
\sum_{i=1}^{l-1} (h_{i,i+1}^1+h_{i,i+1}^2+h_{i,i+1}^3)(a,b,\sigma).$$
\end{lem}

We first prove Lemma {\ref{sum}}.

As an immediate consequence of the definitions (\ref{definehk}) we have:
\begin{eqnarray*}
 & &h_{i,i}^1-h_{i,i+1}^1 =  \sum_{k\leq s_i}
(a_{i,k}+b_{i,k}+1)\ve_{k,(i,i+1]}
\\ & &h_{i,i}^2- h_{i-1,i}^2=\sum_{k\leq s_i}
(a_{i,k}+b_{i,k})\ve_{(i-1,i],k}
\\ & &h_{i,i}^3 - h_{i-1,i}^3 = \sum_{s_{i-1}<k\leq s_i} b_{i,k}
\end{eqnarray*}
and since $h_{l,l+1}^1=h_{0,1}^2=h_{0,1}^3=0$, Lemma \ref{sum}  follows.
\epf 

It only remains to prove Lemma \ref{kk}. We have 
\begin{eqnarray*}
\lefteqn{\Hom(\A_i,\B_j) = \bigoplus_{k,m} \Hom(\A_{i,k},\B_{j,m}) =}\\
& & \bigoplus_{k\leq s_i} 
\bigoplus_{m>s_j}
\Hom(\oh_{\sigma(k)}(-a_{i,k}-b_{i,k}), \oh_{\sigma(m)}))\\
& &\oplus 
\bigoplus_{k\leq s_{i},m\leq s_j}
\Hom(\oh_{\sigma(k)}(-a_{i,k}-b_{i,k}),
   \oh_{\sigma(m)}/ \oh_{\sigma(m)}(-a_{j,m}-b_{j,m}))
\end{eqnarray*}
We have three situations to consider:

\noindent 1. $m>s_j$

$\Hom(\oh_{\sigma(k)}(-a_{i,k}-b_{i,k}), \oh_{\sigma(m)})$ is of rank
$a_{i,k}+b_{i,k}+1$, with weights the same sign as $(\sigma(m)-\sigma(k))$
for $N$ large enough.  
The number of $m>i$ with $\sigma(m)>\sigma(k)$ is 
$\sum_{i<m\leq n} \epsilon_{k,m}$.  In the notation of (\ref{sigma}), since
$s_{l+1}=n$, this number is $(a_{i,k}+b_{i,k}+1)\ve_{k,(j,l+1]}$.
This gives the term $h^1_{i,j}(a,b,\sigma)$ 
in Lemma \ref{kk}.

\vspace{2 ex}\noindent 2.  $k\neq m$

$\Hom(\oh_{\sigma(k)}(-a_{i,k}-b_{i,k}),
   \oh_{\sigma(m)}/ \oh_{\sigma(m)}(-a_{j,m}-b_{j,m}))$ is of rank
$a_{j,m}+b_{j,m}$, with weights the same sign as $\sigma(m)-\sigma(k)$.
Thus the positive weight contribution is
$\sum_{k\leq s_i} (a_{j,m}+b_{j,m})\epsilon_{k,m}=
(a_{j,m}+b_{j,m})\epsilon_{(0,i],m}$, which gives the term
$h^2_{i,j}$ in Lemma \ref{kk}.

\vspace{2 ex}
\noindent 3.  $k=m$

$\Hom(\A_{i,k},\B_{j,k})$ sits inside the long exact sequence induced by
the short exact sequence (\ref{ses}):
\begin{eqnarray}
\label{les}
0 &\rarr& \Hom(\A_{i,k},\A_{j,k})\rarr \Hom(\A_{i,k},\oh_{\sigma(k)})\\
\nonumber & &\mbox{   }\rarr \Hom(\A_{i,k},\B_{j,k})
 \rarr H^1(\A_{i,k}^*\otimes\A_{j,k}) \rarr 0
\end{eqnarray}
We have two cases:  

\vspace{2 ex}
(a)  $\A_{i,k}\hookrightarrow \A_{j,k}$.  Here, equation (\ref{les}) becomes
$$0 \rarr \Hom(\A_{i,k},\A_{j,k}) \rarr
\Hom(\A_{i,k},\oh_{\sigma(k)}) 
\rarr \Hom(\A_{i,k},\B_{i,k}) \rarr 0. $$
$\Gamma$ acts on $\Hom(\A_{i,k},\A_{j,k})$ by $t\cdot X^rY^s =
X^{r-(a_{i,k}-a_{j,k})}Y^{-s+(b_{i,k}-b_{j,k})}$ and on 
$\Hom(\A_{i,k},\oh_{\sigma(k)})$
 by $t\cdot X^rY^s = X^{r-a_{i,k}}Y^{-s+b_{i,k}}$.  
Thus, the positive part of this piece of $T_{\R}(r(a,b,\mathbf{\sigma}))$
has dimension $b_{j,k}$.  This gives the $h^3_{i,j}$ term.

\vspace{2 ex}
(b)  $\A_{i,k}\not\hookrightarrow \A_{j,k}$. From (\ref{les}) we get
$$0 \rarr \Hom(\A_{i,k},\oh_{\sigma(k)})
\rarr \Hom(\A_{i,k},\B_{j,k}) \rarr H^1(\A_{i,k}^*\otimes\A_{j,k}) \rarr 0.$$
We have $H^1(\A_{i,k}^*\otimes\A_{j,k})=H^1(\oh((a_{i,k}-a_{j,k})+
(b_{i,k}-b_{j,k}))).$ 
By Serre duality and the same arguments as in the previous case, the
positive contribution of $\Hom(\A_{i,k},\B_{j,k})$ is 
$(b_{j,k}-b_{i,k})+b_{i,k}=b_{j,k}$.  This gives the term $h^3_{i,j}$.

Therefore we have shown that
$\dim_\com (T_{\R}(\A))^+ = h(a,b,\sigma)$, so that the proposition is proved.
\epf

\section{Poincar\'{e} polynomials}
\label{poincare}
We  use the implicit formula for the Betti numbers proved
in Proposition \ref{bpf}.
Rewrite (\ref{hpf}) as:
\begin{eqnarray}
\label{rehpf} 
 \lefteqn{h(a,b,\sigma) = 
\sum_{i=1}^l \sum_{k\leq s_i} \ve_{k(i,i+1]}
+ \sum_{i=1}^l \sum_{s_{i-1}<k\leq s_{i}} b_{i,k}
 }\\
\nonumber & &\mbox{}+ 
\sum_{i=1}^l \sum_{k\leq s_i} 
(a_{i,k}+b_{i,k})(\ve_{k(i,i+1]}+ \ve_{(i-1,i]k})
\end{eqnarray}

Recall the sequences of independent nonnegative integers $\ta$ and $\tb$
introduced in the proof of Theorem {\ref{eulerthm}},
given by $\ta_{i,j}=a_{i,j}-a_{i+1,j},\tb_{i,j}=b_{i,j}-b_{i+1,j}$.
We see that $a_{i,k} = \sum_{j\geq i} \ta_{j,k}$ and
$b_{i,k} = \sum_{j\geq i} \tb_{j,k}$.
Changing to the variables  $(\ta,\tb,\sigma)$, the middle sum of
(\ref{rehpf}) becomes
\begin{eqnarray*}
\sum_{i\leq j\leq l} \sum_{k\leq s_i} 
(\ta_{j,k}+\tb_{j,k})(\ve_{k(i,i+1]}+ \ve_{(i-1,i]k})\\
= \sum_{i\leq j\leq l} \sum_{s_{i-1}<k\leq s_i} 
(\ta_{j,k}+\tb_{j,k})(\ve_{k(i,j+1]}+ \ve_{(i-1,j]k})
\end{eqnarray*}

Therefore, we can simplify our expressions to get:
\begin{eqnarray}
\label{Hw}
\lefteqn{ H(\ta,\tb,\sigma) :=  h(a,b,\sigma) = 
\sum_{i=1}^l \sum_{k\leq s_i} \ve_{k(i,i+1]}
+ \sum_{j=1}^l \sum_{k\leq s_j}  \tb_{j,k} }\\ \nonumber
& &\mbox{} +\sum_{i\leq j\leq l} \sum_{s_{i-1}<k\leq s_i}
(\ta_{j,k}+\tb_{j,k})(s_j-s_i + \ve_{k(j,j+1]}+\ve_{(i-1,i]k})
\end{eqnarray}

By the definition of $H$, we have shown
\begin{prop}
\label{Hpf}
$$b_{2M}(\Rpf) =\mbox{rank }A_M(\R) = \#(H^{-1}(M)).$$
\end{prop}

For $w\in S$, define $P(w)$ to be the elements 
$(\ta,\tb,\sigma)\in P$ satisfying
$\sigma = w$.  Let $H_w$ denote the restriction of $H$ to
$F(w)$. Then $\#(H_w^{-1}(M))$ counts the
number of sequences of natural numbers $(\ta,\tb)$ satisfying
the  relation $H_w(\ta,\tb) = M$ given by (\ref{Hw}) as well as
the relations in (\ref{rel}).  Since all of these relations are linear
in the variables
$\ta_{k,j}$ and $\tb_{k,j}$, by the same reasoning as in the proof 
of Theorem \ref{eulerthm}, we have the following generating function:

\begin{eqnarray}
\label{pfw}
\lefteqn{ \sum_{M,d_1,...d_l} \#(H_w^{-1}(M))t_1^{d_1}\cdots t_l^{d_l}z^M = }\\
& &\nonumber z^{\sum_{i\leq l} \sum_{k\leq s_i} \ve^w_{k(i,i+1]}}
\prod_{1\leq i \leq j \leq l} \prod_{s_{i-1}<k\leq s_i}
\frac{1}{f^{i,j}_{\rho_{i,j,k}}}
\frac{1}{f^{i,j}_{\rho_{i,j,k}+1}}
\end{eqnarray}
where we have defined $f^{i,j}_k= 1-t_i...t_jz^k$ and
$\rho_{i,j,k} = s_j-s_i + \ve^w_{k(j,j+1]}+k-s_{i-1}-1.$

We are now ready to prove Theorem 1.

By the definitions of $H$ and $H_w$,
we have $$\#(H^{-1}(M)) = \sum_{w\in S} \#(H_w^{-1}(M))$$ and hence
$b_{2M}(\Rpf)= \sum_w \#(H_w^{-1}(M))$.  Therefore, by (\ref{pfw}),
it suffices to prove the following (purely combinatorial) result: 
\begin{prop}
\label{comb}
\begin{eqnarray*}
\sum_{w\in S(n;s_1,...,s_l)} 
z^{\sum_{i\leq l} \sum_{k\leq s_i} \ve^w_{k(i,i+1]}}
\prod_{1\leq i \leq j \leq l} \prod_{s_{i-1}<k\leq s_i}
\frac{1}{f^{i,j}_{\rho_{i,j,k}}}
\frac{1}{f^{i,j}_{\rho_{i,j,k}+1}}\\
=
\left(\frac{\prod_{i=1}^n (1-z^i)}{\prod_{j=1}^{l+1} 
\prod_{i=1}^{s_j-s_{j-1}} (1-z^i)} \right)
\prod_{1\leq i \leq j \leq l} \prod_{s_{i-1}<k\leq s_i}
\frac{1}{f^{i,j}_{s_j-k }}
\frac{1}{f^{i,j}_{s_{j+1}-k+1}}
\end{eqnarray*}
\end{prop}

We use induction on $n$ to prove the proposition.  

For $n=1$, there are two cases:
\begin{enumerate}
\item $s_1 = 0$.  Here, $S(1;0)=S_1={id}$.  Both sides of the equation are
equal to $1$, so that the proposition holds.

\item $s_1 = 1$.  We have $S(1;1)=S_1={id}$.  Both sides of the equation
are equal to $\frac{1}{(1-z)(1-tz)}$, so again the proposition holds.
\end{enumerate}

The strategy is to break up $S=S(n;\s)$ into $l+1$ different permutation
groups upon which we can use the inductive hypothesis.

For $1\leq m \leq l+1$, let $S(m)$ denote the subset of $S$ consisting
of permutations $w$ such that $w(s_{m-1}+1)=1$.  It is clear that 
$S = \bigcup_m S(m)$ is a disjoint union.  Any $w\in S(m)$ satisfies:
$\epsilon^w_{s_{m-1}+1,j} = 1$ for $j\neq s_{m-1}+1$ and 
$\epsilon^w_{i,s_{m-1}}=0$ for all $i$.

For $w\in S(m)$, let $w'\in S':=S(n-1;s'_1,...,s'_l)$ be defined by:
\begin{enumerate}
\item
$w'(i) = w(i)-1$ for $i\leq s_{m-1}$
\item
$w'(j) = w(j+1)-1$ for $s_{m-1}+1<j$.   
\end{enumerate}

Let $\epsilon'_{i,j}=\epsilon^{w'}_{i,j}$ and define $\ve'_{i(k,k']}$ and
$\rho'_{i,j,k}$ accordingly.
Note that $s'_i = s_i$ for  $i\leq m-1, s'_j = s_j-1$ for $m \leq j$ and
\begin{equation}
\label{epsilonsum}
\sum_{i\leq l} \sum_{k\leq s_i} \ve _{k(i,i+1]}=
(n-s_m) + \sum_{i\leq l} \sum_{k\leq s_i, k\neq s_{m-1}+1} \ve_{k(i,i+1]}.
\end{equation}

We split the left hand sum over $S$ into sums over $S(m)$ for $1\leq m
\leq l+1$.  For each $S(m)$, we can factor out the parts of the
product on the left hand side  where  $i=m$ and $k=s_{m-1}+1$.  Using
(\ref{epsilonsum}), this gives:

\begin{eqnarray*}
\lefteqn{ \sum_{w\in S} z^{\sum_{i\leq l} \sum_{k\leq s_i} \ve_{k(i,i+1]}}
\prod_{1\leq i \leq j \leq l} \prod_{s_{i-1}<k\leq s_i}
\frac{1}{f^{i,j}_{ \rho_{i,j,k}}}
\frac{1}{f^{i,j}_{ \rho_{i,j,k}+1}}  }\\
&=&\sum_{m=1}^{l+1} z^{n-s_m} 
\prod_{m\leq j}
\frac{1}{f^{m,j}_{s_{j+1}-s_m}}
\frac{1}{f^{m,j}_{s_{j+1}-s_m+1}} \\
& & \left( \mbox{}\sum_{w\in S(m)}
z^{\sum_{i\leq l} \sum_{k\leq s'_i} \ve'_{k(i,i+1]}}
\prod_{i\leq m\leq j} \prod_{s'_{i-1}<k\leq s'_i}
\frac{1}{f^{i,j}_{ \rho'_{i,j,k}+1}}
\frac{1}{f^{i,j}_{ \rho'_{i,j,k}+2}} \right. \\
& & \left.
\prod_{\tiny{ \begin{array}{c} i\leq j \\ m>i \mbox{ or } j<m \end{array} }}
\prod_{s'_{i-1}<k\leq s'_i}
\frac{1}{f^{i,j}_{ \rho'_{i,j,k}}}
\frac{1}{f^{i,j}_{ \rho'_{i,j,k}+1}} \right).
\end{eqnarray*}

By induction. we may apply the result to 
$S(m)\eqq S(n-1;s'_1,\ldots,s'_{l+1})$
to the quantity in parentheses.  In particular, by what follows from the 
proof of Proposition \ref{comb} for $S(m)$, the sum becomes:

\begin{eqnarray*}
\lefteqn{ \sum_{m=1}^{l+1} z^{n-s_m} 
\prod_{m\leq j}
\frac{1}{f^{m,j}_{s_{j+1}-s_m}}
\frac{1}{f^{m,j}_{s_{j+1}-s_m+1}} }\\
& &\left( \frac{\prod_{i=1}^{n-1} (1-z^i)}{\prod_{j=1}^{l+1} 
\prod_{i=1}^{s'_j-s'_{j-1}} (1-z^i) }
 \prod_{i\leq m \leq j}
\prod_{s'_{i-1}<k\leq s'_i}
\frac{1}{f^{i,j}_{s'_j-k}}
\frac{1}{f^{i,j}_{s'_{j+1}-k+1}} \right.  \\
& &\left. 
\prod_{\tiny{ \begin{array}{c} i\leq j \\ m>i \mbox{ or } j<m \end{array} }}
\prod_{s'_{i-1}<k\leq s'_i}
\frac{1}{f^{i,j}_{s'_j +1-k}}
\frac{1}{f^{i,j}_{s'_{j+1}+1-k+1}} \right).
\end{eqnarray*}

Since we have:
$$ \frac{\prod_{i=1}^n (1-z^i)}{\prod_{j=1}^{l+1} 
\prod_{i=1}^{s_j-s_{j-1}} (1-z^i)} =
\frac{\prod_{i=1}^{n-1} (1-z^i)}{\prod_{j=1}^{l+1} 
\prod_{i=1}^{s'_j-s'_{j-1}} (1-z^i)} \frac{1-z^n}{1-z^{s_m-s_{m-1}}},$$
we can write this sum explicity as:
\begin{eqnarray*}
& &\sum_{m=1}^{l+1} z^{n-s_m}  \frac{\prod_{i=1}^n (1-z^i)}{\prod_{j=1}^{l+1} 
\prod_{i=1}^{s_j-s_{j-1}} (1-z^i)}
\cdot \frac{(1-z^{s_m-s_{m-1}})}{(1-z^n)}
\\ & &\prod_{m\leq j} 
\frac{1}{f^{m,j}_{s_{j+1}-s_m}}
\frac{1}{f^{m,j}_{s_{j+1}-s_m+1}}
\cdot \prod_{m\leq j} \prod_{s_{m-1}<k\leq s_m-1} 
\frac{1}{f^{m,j}_{s_j-k}}\frac{1}{f^{m,j}_{s_{j+1}+1-k}} \\ & &
\prod_{i\leq m-1}
\prod_{s_{i-1}<k\leq s_i} 
\frac{1}{f^{i,m-1}_{s_{m-1}-k}}\frac{1}{f^{i,m-1}_{s_{m}-1-k}} 
\prod_{\tiny{ \begin{array}{c} i\leq j \\ i\neq m \\ j \neq m-1 
\end{array} }}
\prod_{s_{i-1}<k\leq s_i} 
\frac{1}{f^{i,j}_{s_j -k}}\frac{1}{f^{i,j}_{s_{j+1}+1-k}}
\end{eqnarray*}

By clearing denominators, the proposition follows once we have proven:

\begin{lem}
\label{pfpoly}
\begin{eqnarray*}
\lefteqn{ (1-z^n) \prod_{1\leq i \leq j \leq l} f^{i,j}_{s_{j+1}-s_i} =} \\
& &
\sum_{m=1}^{l+1} z^{n-s_m} (1-z^{s_m-s_{m-1}})
\prod_{i \leq m-1} f^{i,m-1}_{s_m-s_{i-1}}
\prod_{m \leq j} f^{m,j}_{s_j-s_m}
\prod_{\tiny{ \begin{array}{c} i\leq j \\ i\neq m \\ j \neq m-1 
\end{array} }} f^{i,j}_{s_{j+1}-s_i}
\end{eqnarray*}
\end{lem}

\bpf
First, we change our notation so that we  work with independent
variables $x_i$ where we define 
$x_r = t_1...t_r$, with $x_0=1$.  Then for any $i \leq j$,
we have $t_i...t_j = x_j/x_{i-1}$.  Therefore we have 
$$f^{i,j}_k = 1- t_i...t_jz^k = \frac{x_{i-1}-x_jz^k}{x_{i-1}}.$$

Substituting into both sides of Lemma \ref{pfpoly} and multiplying
through by $x_0^{l}x_1^{l-1}...x_l^{0}$ 
it suffices to prove the following
polynomial identity in the ring $\com[x_0,...x_n,z]$,
where we define $e^{i,j}_{k} = x_i-x_jz^k$:
\begin{equation}
\label{xpf}
(1-z^n) \prod_{1\leq i \leq j\leq l} e^{i-1,j}_{s_{j+1}-s_i } =
\end{equation}
$$ \sum_{m=1}^{l+1} z^{n-s_m}(1-z^{s_m-s_{m-1}})
 \prod_{i \leq m-1} e^{i-1,m-1}_{s_m-s_{i-1}}
\prod_{m \leq j} e^{m-1,j}_{s_j-s_m}
\prod_{\tiny{ \begin{array}{c} i\leq j \\ i\neq m \\ j \neq m-1 
\end{array} }} e^{i-1,j}_{s_{j+1}-s_i}.$$

The polynomial ring $\com[x_0,...x_n,z]$ is a unique factorization
domain, and that the left side of (\ref{xpf}) completely factored
except for the term $(1-z^n)$.
Since the degree of $z$ matches on both sides,
as does the term of $t_1^0...t_{l}^0=1$,
namely $(1-z^n)$,
it is enough to show the right side vanishes
with the relation $e^{i-1,j}_{s_{j+1}-s_i}=0$ for each $i\leq j$, i.e.
that $(x_{i-1}-x_jz^{s_{j+1}-s_i})$ is a factor of the right side.

Fix $i\leq j$.  Then all but two summands on the right vanish.  In particular, 
after some cancellations in the final terms of the two products, and
bringing one of the terms to the other side, we
have left to show:
\begin{eqnarray}
\label{pfpolyid}
z^{n-s_i} \prod_{p \leq i-1} e^{p-1,i-1}_{s_i-s_{p-1}}
\prod_{i \leq q} e^{i-1,q}_{s_q-s_i}
\prod_{\tiny{ \begin{array}{c} p\leq j \\ p\neq i
\end{array} }} e^{p-1,j}_{s_{j+1}-s_p}
\prod_{ j+1\leq q} e^{j,q}_{s_{q+1}-s_{j+1}} \\
 \nonumber  =
- z^{n-s_{j+1}} \prod_{p \leq j+1} e^{p-1,j}_{s_{j+1}-s_{p-1}}
\prod_{j+1 \leq q} e^{j,q}_{s_q-s_{j+1}}
\prod_{\tiny{ \begin{array}{c} i\leq q \\ q\neq j 
\end{array} }} e^{i-1,q}_{s_{q+1}-s_i}
\prod_{p\leq i-1 } e^{p-1,i-1}_{s_i-s_p}
\end{eqnarray}
with the relation $e^{i-1,j}_{s_{j+1}-s_i}=0$.

But now, by doing the substitution
$x_{i-1}=x_jz^{s_{j+1}-s_i}$ ,
 we see that we will have proven Lemma \ref{pfpoly} 
once we show that the polynomial identity (\ref{pfpolyid}) holds  
in the (unique factorization domain) $\com[x_0,...x_{i-2},x_i,...x_n,z]$.
We make the following observations about the substitution:

\begin{enumerate}
\item
$e^{p-1,i-1}_{s_i-s_{p-1}} \mapsto e^{p-1,j}_{s_{j+1}-s_{p-1}}$
\item 
$e^{p-1,i-1}_{s_i-s_{p}} \mapsto e^{p-1,j}_{s_{j+1}-s_{p}}$
\item 
$e^{i-1,q}_{s_q-s_i} \mapsto  
\left\{ \begin{array}{ll} z^{s_{j+1}-s_i} e^{j,q}_{s_q-s_{j+1}} 
& \mbox{when } j+1\leq q\\
-z^{s_q-s_i} e^{q,j}_{s_{j+1}-s_q} &  \mbox{when } q \leq j  \end{array} \right. $
\item 
$e^{i-1,q}_{s_{q+1}-s_i} \mapsto
\left\{ \begin{array}{ll} z^{s_{j+1}-s_i} e^{j,q}_{s_{q+1}-s_{j+1}} 
& \mbox{when } j\leq q\\
-z^{s_{q+1}-s_i} e^{q,j}_{s_{j+1}-s_{q+1}} &  \mbox{when } q < j  \end{array} \right. $\end{enumerate}

Substituting into both sides of (\ref{pfpolyid}), and using the
 above properties, we  have two completely
factored polynomials on each side of the identity.
It is easy to check that the degree of $z$ in the two terms match,
that the sign matches, and that the factors of the form $e^{i,j}_k$
match exactly.
\epf

This concludes the proof of Proposition \ref{comb} and Theorem \ref{Rpf}.

For any scheme $X$, setting $z=1$ into
the Poincar\'{e} polynomial $$\poinc(X) = \sum_M (-1)^M b_M(X)z^M$$
gives the Euler characteristic
$\chi(X)$.  Since odd cohomology of $\R$ vanishes,
this substitution into Theorem \ref{Rpf}
provides another proof of Theorem \ref{eulerthm}.

\section{Special cases}
\label{special}

\subsection{Quot scheme}
We can apply Theorem \ref{Rpf} to the ordinary Quot scheme, 
 parametrizing rank $r$ degree $d$  quotients of $V_\Pone$, to get a
generating function for the Poincar\'{e} polynomials of $\Rg$.

\begin{thm}
\begin{eqnarray*}
\lefteqn{\sum_d \poinc(\Rg)t^d = }\\
& & \poinc(\grass^r(n))\cdot 
\prod_{i=1}^{n-r} \left(\frac{1}{1-tz^{i-1}}\right)
\left(\frac{1}{1-tz^{n-i+1}}\right).
\end{eqnarray*}
\end{thm}
\noindent 
where $\poinc(\grass^r(n)$ is the Poincar\'{e} polynomial of $\grass^r(n)$, which
is the following  classical generating 
function for the Betti numbers for the Grassmannian:

\begin{eqnarray*}
\poinc(\grass^r(n)) =
\sum_M b_{2M}(\grass^r(n))z^M =
\frac{\prod_{i=1}^n (1-z^i)} 
{\prod_{i=1}^{n-r} (1-z^i)\prod_{i=1}^r (1-z^i)} 
\end{eqnarray*}

This was the case studied by Str{\o}mme, who found implicit
formulas for the Betti numbers, which are the same up to notation as those 
found in this paper.  Generators and relations of the Chow ring of 
$\R$ are also given  in \cite{s}.  However, this set is far from minimal, and
is not suited to the study of the Chow rings as the degree 
$d$ becomes large.

\subsection{Hyperquot scheme for a complete flag variety}
Consider the 
 space $\Rpf$ where $l=n-1$ and $s_i=i$, i.e. the space $\Rcf$.

Fix $n$.  Then we have the following generating function for the Poincar\'{e}
polynomials of $\Rcf$ where $\bd = (d_1,...,d_{n-1})$.

\begin{thm}
\label{Rcf}
\begin{eqnarray*}
\lefteqn{ \sum_{d_1,...d_{n-1}} \poinc(\Rcf)t_1^{d_1}...t_{n-1}^{d_{n-1}}= }
 \\
& &\poinc(\F(n)) \cdot
\prod_{1\leq i \leq j \leq n-1}
\left(\frac{1}{1-t_i...t_jz^{j-i}}\right)
\left(\frac{1}{1-t_i...t_jz^{j-i+2}}\right).
\end{eqnarray*}

\end{thm}
\noindent
The classical term $\poinc(\F(n))$ is given by: 
$$\poinc(\F(n))=\sum_M b_{2M}(\F(n))z^M = \frac{\prod_{i=1}^n (1-z^i)}
    {\prod_{j=1}^{n} (1-z) }.
$$

\subsection{Examples}

\begin{enumerate}
\item
$\F(1;0) = \grass^1(1)$ is a point and $\mathcal{HQ}_d({\F(1;0)})$ is a point for $d=0$
and empty for $d>0$.
$$\sum_{d}\poinc(\mathcal{HQ}_d(\F(1;0)))t^d = 1$$
which is consistent with the theorem.
\item
 $\F(1;1) = \grass^0(1)$ is a point. $\mathcal{HQ}_d(\F(1;1))$
parametrizes quotients of $\oh_\Pone$ of rank $0$ and degree $d$,
which are all of the form $\oh \rrarr \oh_D$, where  
$D\in \mbox{Sym}^d \Pone \eqq \proj^d$.  Therefore 
$\mathcal{HQ}_d(\F(1;1))  \eqq \proj^d$, giving:
$$\sum_{d,M} b_{2M}(\proj^d)t^dz^M = \sum_d t^d
\left( \sum_{M\leq d} z^M \right) = \frac{1}{(1-t)(1-tz)}.$$
\item 
$\mathcal{HQ}_d(\grass^{n-1}(n))$ 
can be viewed as the space of sheaf injections
$\oh(-d)\hrarr \sumo_{i=1}^{n} \oh$ up to equivalence.  Each inclusion
of sheaves is given by $n$ sections in $H^0(\proj^1,\oh(d))$.  Thus, we
can view any such inclusion as an element of the vector space
$\sumo_{i=1}^{n} H^0(\oh(d))$ of dimension $n(d+1)$.  Two inclusions
are equivalent exactly when they differ by a scalar.  Hence, 
$\mathcal{HQ}_d\eqq
\proj^{n(d+1)-1}$ so that $\poinc(\mathcal{HQ}_d) = \sum_{0\leq M\leq n(d+1)} z^M$.
\end{enumerate}

\end{document}